\documentclass[12pt]{article}
\usepackage{amsfonts, amssymb, amsmath}

\def\R{{I\!\!R}}

\def\CC{{\rm \kern.24em \vrule width.02em height1.4ex
    depth-.05ex \kern-.26em C}}

\def\TagOnRight

\def\AA{{\it I}\hskip-3pt{\tt A}}

\def\QQ{\rlap {\raise 0.4ex \hbox{$\scriptscriptstyle |$}}
  {\hskip -0.1em Q}}

\newcommand{\be}{\begin{equation}}
\newcommand{\ee}{\end{equation}}
\newcommand{\bea}{\begin{eqnarray}}
\newcommand{\eea}{\end{eqnarray}}
\newcommand{\Bea}{\begin{eqnarray*}}
\newcommand{\Eea}{\end{eqnarray*}}

\catcode`\@=11
\def\theequation{\@arabic{\c@equation}}
\catcode`\@=12

\newcommand{\bi}{\begin{itemize}}
\newcommand{\ei}{\end{itemize}}

\newtheorem{Definition}{Definition}[section]
\newtheorem{Theorem}[Definition]{Theorem}

\newtheorem{Proposition}[Definition]{Proposition}

\newtheorem{Remark}[Definition]{Remark}

\renewcommand{\theequation}{ \thesection.\arabic{equation} }

\begin{document}

\title{ Asymptotic distribution of Brownian Excursions
into an Interval.}
\author{B. Rajeev \\
Indian Statistical Institute, 8th Mile,\\
Mysore Road, Bangalore 560 059, India \\
email:brajeev@isibang.ac.in} \maketitle
\begin{abstract} In this paper, following earlier results in \cite{AR}
 we derive the asymptotic distribution, as $t \rightarrow \infty $,
 of the excursion of Brownian motion stradling t,into an interval $(a,b)$ , conditional on
 the event that there is such an excursion.
\end{abstract}
\text{Key words and phrases :} Brownian motion, excursions, last
exit/entrance times, asymptotic distribution of excursions.

\section{Introduction :} In this paper we consider the asymptotic
distribution of the excursions $(\zeta_t)$ of a one dimensional
standard Brownian motion $(W_t)$, straddling the time t, into an
interval $(a,b) $ as $t \rightarrow \infty$. Although, excursions
straddling a given time are well studied in the literature for
general Markov processes (for a sample, see
\cite{M},\cite{P},\cite{RY},\cite{GS}), the study of their
asymptotics as $t \rightarrow \infty$ seems to be new.

To describe our results in more detail let for each $t > 0$,
$\sigma_t,d_t$ be the last entrance before $t$ into $(a,b)$ and
first exit after $t$ from $(a,b)$ respectively for a sample path
$W_\cdot$ such that at time $t, W_t \in (a,b)$. The excursion
straddling time $t$ is the portion of the trajectory $\zeta_t(s) :=
W_{\sigma_t + s\wedge d_t - \sigma_t},  s \geq 0$.  We view the
process $(\zeta_t)$ as a process with values in the space
$C([0,\infty),[a,b])$, the space of continuous functions with values
in $[a,b]$ so that the convergence in question reduces to weak
convergence in this space.The crucial step in the proof is to
express the expected value $E[f(\zeta_t)| W_t \in (a,b)]$, where $
f: C([0,\infty),[a,b]) \rightarrow {\mathbb R}$ a bounded and
continuous function as $E[q(W_{\sigma_t},t-\sigma_t,f)|W_t \in
(a,b)]$ where the kernel $q(x,s,A)$ is a bounded continuous function
of $(x,s), s
> 0, x = a ~{\rm or}~ b $ for a Borel set $A \in C([0,\infty),[a,b])$
and then use the weak convergence of the pair
$(W_{\sigma_t},t-\sigma_t)$. It is known (see \cite{AR}) that the
latter pair converges to $(X,Y)$, the exit place and time
respectively of a Brownian motion started uniformly in the interval
$(a,b)$. Together with the explicit form of $q(x,s,A)$, this gives
the limiting distribution of $\zeta := \lim\limits_{t \rightarrow
\infty}\zeta_t$ as follows : Starting from $X = a~{\rm or}~b$ with
probability $\frac{1}{2}$, the conditional distribution of $\zeta
\in A$ given that the lifetime $Y$ of $\zeta$ is at least $s$, is
given by $q(x,s,A)$. Together with the known distribution of the
lifetime of excursions into $(a,b)$, starting from $a ~{\rm or ~}b$
this gives a complete description of the limiting distribution.

In Section 2, we set up the necessary machinery from excursion
theory. Rather than use the extension of It\^{o}'s excursion theory
, due to B.Maisonneuve (\cite{M}) , in cases where the boundary of
the excursion set involves more than one point, we give a more
intutive,`bare hands' construction, using It\^{o}'s original result
(as presented in \cite{RY})(Theorem 2.3). Although we are dealing
with a very basic example as far as excursion theory is concerned,
our approach maybe of interest in constructing non trivial `exit
systems' starting with excursions from a single point. See Remark
2.4 for the connection with the exit system formalism of \cite{M}.
As mentioned above, the crucial point in the proof is to express the
expected value of functionals of the excursions into $(a,b)$ at time
t, in terms of the kernel $q(x,s,A)$ describing the conditional
excursion measures and then the continuity in $(x,s)$ of these
kernels. This is done via a `conditional excursion formula' for the
excursion straddling $t$, proved in Section 3. Such formulas are
well known for general Markov processes (see \cite{GS},\cite{M}) but
at specific time instances and not for the excursion process as a
whole, as in \cite{RY} and as required in our situation.

In section 4, we describe the limiting distribution of $\zeta_t$ and
prove the convergence to this distribution. As mentioned above the
proof involves the convergence, as $t \rightarrow \infty$ of the
pair of variables $(W_{\sigma_t},t-\sigma_t)$ to the pair $(X,Y)$
described above. This latter result is a consequence of a more
general result, proved in \cite{AR}, about the convergence as $t
\rightarrow \infty$ of the time reversal $(W_{t-s.}, 0 \leq s \leq
t)$ of Brownian motion $W$, from a point $t$ with $W_t \in (a,b)$.
Finally, we give an application to the evaluation of asymptotic
distribution of functionals of the excursion straddling $t$ as $t
\rightarrow \infty$.

\section{Excursions into an interval :}
Let $S$ be an open subset of ${\mathbb R}$. In this paper, $S$ will
be either $(0,\infty)\bigcup(-\infty,0)$ or the  finite interval
$(a,b)$ or the set $(-\infty,a)\bigcup (a,b)\bigcup (b,\infty)$.
$\partial S$ will denote the boundary of $S$. We denote the closure
of a set S by $\bar{S}$. Let $C([0,\infty),\bar{S})$ denote the
space of continuous functions from $[0,\infty)$ into $\bar{S}$,
equipped with the topology of uniform convergence on compact subsets
of $[0,\infty)$. We will denote by ${U} = U(S)$ the space of
excursions from $\partial S$ into $S$. In other words, \Bea U &:=&
\{u \in C([0,\infty),\bar{S}): {\rm u (0) \in
\partial S , and~ \exists R(u)>0 ~ such ~that }\\&& u(t) \in S , 0
<t< R(u) {\rm ~and} ~u(t) \in
\partial S ~\forall~ t \geq R(u)\}.\Eea

Let ${\cal U}$  be the trace of the Borel $\sigma$-field ${\cal B}$
of $C([0,\infty),\bar{S})$ on $U$ i.e. ${\cal U}=U \cap {\cal B}$.
Let $U_\delta := U  \cup \{\delta\}$, with $\delta$ attached as an
isolated point. Let ${\cal U}_{\delta}$ be the sigma field generated
by ${\cal U}$ and $\{\delta\}$.

Let $(\Omega, {\cal F},P)$ be a probability space and  $({\cal F}_t)
$ a filtration on it satisfying the usual conditions i.e. ${\cal
F}_0$ contains $P$ null sets and the filtration is right continuous.
Let $(W_t)$ be a 1-dimensional standard ${\cal F}_t$- Brownian
motion. ${\cal P}$ will denote the previsible $\sigma$-field over
$({\cal F}_t)$. Let $\left(L^x_t\right)_{t\geq0}$ denote the local
time process of $(W_t)_{t\geq 0}$ at $x\in \R$. For a continuous ,
non decreasing , ${\cal F}_t$-adapted process $(L_t)$ , we define
its right continuous inverse $(\tau_t)$ in the usual way by
\[
\tau_t :=\inf \{s>0 : L_s > t\}.
\]
We note that $(\tau_t)$ is a non decreasing, right continuous
process such that for each $t, \tau_t$ is an ${\cal F}_t$-stopping
time. Let $Z$ be the random closed set $Z := \{t : W_t \in S^c\}$.
We will assume that a.s. (P), $L_t = \int_0^t I_Z(s)~dL_s$. For
$t>0$ let $e_t : \Omega \rightarrow U_{\delta}$ be defined by
\begin{eqnarray*}
e_t (\omega )(s) &=& W_{\tau_{t-}+s\wedge{(\tau_t -
\tau_{t-})}}(\omega),~~
s \geq 0,~~ \tau_t(\omega) -\tau_{t-}(\omega)\neq 0 ,\\
&& \qquad \delta \qquad  ~~~~~~~~~~~~~~~~~~~~~~~\tau_t(\omega)
-\tau_{t-}(\omega)=0.
\end{eqnarray*}

Then $e(.)$ defines a  sigma finite point process on $\Omega$ with
state space ${U}_{\delta}$ and time domain $D(\omega) =\{t>0 :
\tau_{t-} (\omega) \neq \tau_t (\omega)\}$. For $\Gamma \in {\cal
U}_{\delta}$ we define
\[
N_t^{\Gamma} (\omega) := \# \{ s \leq t :s \in D (\omega), e_s
(\omega) \in \Gamma\}.
\]
With $Z$ as defined above we can write $Z(\omega)^c =
\bigcup_{i=1}^{\infty}(\alpha_i (\omega), \beta_i (\omega))$. Let
$G(\omega):=\{\alpha_i (\omega) : i= 1, \ldots \}$ be the set of
left end points of the excursion intervals , $d_t (\omega) :=\inf
\{s>t:W_s (\omega)\notin S \}$ and $\sigma_t = \sup\{s \leq t : W_s
\notin S\}$. We use the convention that $\inf\{\phi\} = \infty$ and
$\sup\{\phi\} = 0.$ For $\omega \in \Omega ,$ and $t
> 0,$ we define the portion of the excursion straddling
$t$~, from $t$ onwards viz. $i_t(\omega)$ as follows : \bea i_t
(\omega)(s) :=W_{t+s\wedge (d_t(\omega)-t})(\omega),~~ 0 \leq s \leq
d_t(\omega)-t. \eea Note that if $t\in G(\omega)$, then $i_t
(\omega) \in U(S)$.

Now we specialise to the case $S = (0,\infty)\bigcup(-\infty,0)$ and
recall a few well known facts. We take $L_t := L_t^0 ,$ the local
time at $0$ of Brownian motion. Then it is well known that $e$ is a
$\sigma$- finite $({\cal F}_{\tau_t})$ Poisson point process with
characteristic measure $n(\Gamma), \Gamma \in {\cal U}_{\delta}$
given by
\[
n (\Gamma) := \frac{1}{t} \, E \left(N^\Gamma_t\right).
\]
Note that $n(\{\delta\})= 0$.
\begin{Remark}
Let $t >0. $ Then it is known (see \cite{RY},Propn.2.8,Chap.XII)
that $n\{u \in U : R(u) > t \} = (\frac{2}{\pi t})^{\frac{1}{2}}$.
In particular $N^\Gamma_t < \infty$ a.s. with \\$\Gamma = \{u : R(u)
> t\}.$
\end{Remark}
We also note that  if $\Gamma $ is such that $EN_t^\Gamma< \infty$
then $n (\Gamma)<\infty$ and \\ $\{N_t^\Gamma-t n (\Gamma) ; t \geq
0 \}$ is an ${\cal F}_{\tau_t}$-martingale.

Suppose $H: [0,\infty) \times \Omega \times U_{\delta} \rightarrow
[0,\infty]$, with $H(t,\omega,\delta)\equiv 0$ for all $(t,\omega)$.
Suppose $H$ is ${\cal P} \otimes {\cal U}_{\delta}$ measurable.
  We have the following well known Theorem of K.It\^{o}(\cite{I}).

\begin{Theorem}
Let $S = (-\infty,0)\bigcup(0,\infty), L_t := L_t^0,$ and $H,G,
i_t,\tau_t$ and $n$ be as above. Then, \bea E\left[\sum\limits_{t
\in G(\omega)} H(t,\omega,i_t(\omega))\right] &=&
E\left[ \sum\limits_{s \in D(\omega)} H(\tau_{s-} (\omega), \omega,e_s(\omega))\right] \nonumber \\
&=& E \int\limits_0^\infty ds \int\limits_{U_{\delta}} H(\tau_s (\omega),\omega,u) n(du)\nonumber \\
&=& E \int\limits_0^\infty dL_s (\omega) \int\limits_{U_{\delta}}
H(s,\omega,u) n (du). \eea
\end{Theorem}

{\bf Proof:} The proof follows as in Proposition 2.6, using
Proposition 1.10  and Theorem 2.4 of Chapter XII of \cite{RY}(see
also \cite{B},Chap.III, Thm.3.18 and Thm 3.24) . The only difference
in our case is that ${\cal F}_t$ is not the canonical filtration.
\hfill{$\Box$}

 ~We now take $S = (a,b), -\infty < a < b < \infty $ and
 $L_t := L_t^a +L_t^b,~~ t\geq
0.$ Let  $(\tau_t)_{t\geq 0},  D $ be as defined above . We note
that for $t \in D(\omega)$ the excursion $e_t(\omega)$ defined
earlier , may not be in $U$ since $D(\omega)$ includes excursions
below $a$ and above $b$. Consequently $G(\omega) \subseteq
D(\omega)$. Since we are only interested in excursions into $(a,b)$
we proceed as follows : For $t \in D(\omega)$, define $e_t'
(\omega)\in U$ as
\begin{eqnarray*}
e_t'(\omega) (s) & := & W_{\tau_{t-}(\omega)+s}(\omega) \qquad 0 \leq s < \tau_{t}(\omega) -\tau_{t-}(\omega)\\
&& W_{\tau_t(\omega)}(\omega)~~~~~ \qquad s \geq \tau_t (\omega)
-\tau_{t-} (\omega).
\end{eqnarray*}
If $t \not\in D(\omega)$, put $e_t'(\omega)\equiv \delta$. Define
$$e_t(\omega) := I_{U}(e_t'(\omega))e_t'(\omega)+ \delta
I_{U^c}(e_t'(\omega)).$$ We shall take this as the definition of our
excursion process when $S = (a,b).$

We define the kernel $n(x,\cdot) , x \in {\mathbb R}$ on $U=U(S) , S
= (a,b)$ as follows: We denote by  $U_0 := U(S),{\rm when}~ S =
(-\infty,0)\bigcup(0,\infty)$, the space of excursions from 0. Let
$U_0^+ :=\{u \in U_0, u(t)
>0, 0 < t<R(u)\}$. $U^-_0$ is defined similarly. $U_0 =U_0^+
\cup U^-_0$ and $U_{0,\delta} := U_0 \bigcup \{\delta \}$. The sigma
field ${\cal U}_{0,\delta}$ is the sigma field generated by ${\cal
U}_0$ and $\{\delta\}$, where ${\cal U}_0 = U_0 \cap {\cal B}$ and
${\cal B}$ the Borel sigma field of $C([0,\infty),{\mathbb R})$. For
$c \in {\mathbb R}$, let $T_c :=\inf \{s>0:W_s =c\}$. For $u \in
U_0$, let $u^c$ denote the path in $C([0,\infty), {\mathbb R})$
which is given by $u(\cdot)$ stopped when it reaches level $c$ viz.
$u(T_c\wedge \cdot)$. Define maps $ \lambda^a :U_{0,\delta}
\rightarrow U_{\delta} ,~ \lambda^b :U_{0,\delta} \rightarrow
U_{\delta }$ as follows: For $u \in U^+_0$, $\lambda^a(u)
:=a+u^{b-a}$ and for $u \in U^-_0,~ \lambda^b(u):=b+u^{a-b}$. We
extend $\lambda^a,\lambda^b$ to the whole of $U_{0,\delta}$ by
setting $\lambda^a = \delta~ (= \lambda^b) $ on $U^-_0 \cup
\{\delta\}$ (respectively $U_0^+ \cup \{\delta\}$).

Let $n$ be the It\^{o} excursion measure and $ n^{+}:=n \mid_{U_0^+}
, n^{-}:=n \mid_{U_0^-}$, the restrictions of $n$ to $U_0^+$ and
$U^-_0$ respectively.

Define $n_a := n^+ \circ (\lambda^a)^{-1}$; $n_b := n^- \circ
(\lambda^{b})^{-1}$. For $\Gamma\in {\cal U}, x \in {\mathbb R}$,
define
$$n(x,\Gamma):= 1_{\{a\}} (x) n_a (\Gamma )+1_{\{b\}} (x)
n_b(\Gamma ).$$ We extend $n(x, .)$ to the whole of ${\cal
U}_{\delta}$ by setting $n(x,\{\delta\})=0 $ for every $x$. Let
$H:[0,\infty) \times \Omega \times U_\delta \rightarrow [0,\infty]$
be ${\cal P} \otimes {\cal U}_\delta$ measurable and such that
$H(t,\omega,\delta) = 0 $ for all $(t,\omega)$.

\begin{Theorem}
Let $S = (a,b)$. Let, $G,H, i_t ,\tau_t,L_t$ and $n(x,\cdot)$ be as
above , corresponding to the interval $(a,b)$. Then,  \bea E\left[
\sum\limits_{t\in G(\omega)} H(t,\omega,i_t (\omega))\right]
&=& E\left[ \sum\limits_{s\in D(\omega)} H(\tau_{s-}(\omega), \omega,e_s (\omega))\right] \nonumber \\
&=&  E\int\limits_0^\infty dL_s(\omega) \int\limits_{U_\delta}
H(s,\omega,u) n(W_s, du) \nonumber \\
&=& E \int\limits_0^\infty ds \int\limits_{U_\delta}
H(\tau_s,\omega,u) n(W_{\tau_s},du). \eea
\end{Theorem}

{\bf Proof:} The first equality follows from the inclusion
$G(\omega) \subseteq D(\omega)$ and the fact that for $s \in
D(\omega)- G(\omega), H(\tau_{s-},\omega,e_s(\omega))= 0$. The third
equality follows by time change . Thus it suffices to prove the 2nd
equality in the statement. Let $D^+(\omega):=\{t\in
G(\omega):W_{\tau_{t-}}=a\}$ and $D^-(\omega):=\{t\in
G(\omega):W_{\tau_{t-}}=b\}$. Then
\begin{eqnarray*}
E \sum\limits_{t \in D} H(\tau_{t-},\omega,e_t) &=&
E\sum\limits_{t \in D^+} H(\tau_{t-},\omega,e_t)+E\sum\limits_{t\in D^-} H(\tau_{t-},\omega,e_t)\\
&=:& S_1 +S_2.
\end{eqnarray*}
To analyse $S_1$ introduce the standard Brownian motion $(\tilde
W_t)$ where
\[
\tilde W_t := (W_{t+T_a} -W_{T_a})I_{\{T_a < \infty \}} =
(W_{t+T_a}-a)I_{\{T_a < \infty \}}~.
\]
 We will denote the excursions from 0 of $(\tilde W_t)$ with a
tilde. Thus, $\tilde L_t^0$ is the local time at 0 of $(\tilde W_t)$
with right continuous inverse $(\tilde \tau_t^0)$; $ \tilde e_t^0$
is the excursion process for $t\in \tilde D^0 =\{s: \tilde
\tau_{s-}^0 \neq \tilde \tau_s^0\}$. Let $\tilde D^{0,+}=\{s\in
\tilde D^0 : \tilde{e}_s^0 \in U^{+}_0 \}$. Then note that, almost
surely, there is a 1-1 correspondence between $t \in D^+(\omega)$
and $s\in \tilde D^{0,+}$ in the sense that $\tau_{t-}(\omega)
=\tilde \tau_{s-}^0 (\omega)+T_a$. This follows from two facts :
Firstly the local time at $0$ for $(\tilde W_t)$ at time $t$ is
precisely the local time at $a$ for $(W_t)$ at time $t + T_a$ on
$T_a < \infty.$ Secondly, the positive excursions of $(\tilde W_t)$
from $0$ until the hitting time of $b-a$ are exactly the excursions
of $(W_t)$ from $a$ until the hitting time of $b.$ Further for such
$t$ and $s$, $e_t (\omega) =\tilde e_s^{0,b-a}+a=\lambda^a (\tilde
e_s^0)$. Hence
\begin{eqnarray*}
S_1 &=& E\sum\limits_{s\in \tilde D^{0,+}} H\left( \tilde
\tau^0_{s-}
 +T_a ,\omega,\lambda^a (\tilde e_s^0)\right)I_{\{T_a < \infty\}}\\
&=& E\sum\limits_{s \in \tilde D^0} \tilde H (\tilde \tau_{s-}^0,
\omega,\tilde e_s^0)
\end{eqnarray*}
where $\tilde H (t,\omega,u) := H(t+T_a,\omega,\lambda^a(u))I_{\{T_a
< \infty\} } $ for $t \geq 0, \omega \in \Omega, u \in
U_{0,\delta}$. Note that $\tilde H (t,\omega,u)$ is $\tilde {\cal
P}\otimes {\cal U}_{0,\delta}$ measurable , where $\tilde{\cal P}$
is the previsible sigma field with respect to the filtration $({\cal
F}_{t+T_a})$. Since $(\tilde{W}_t)$ is an $({\cal F}_{t+T_a})$
Brownian motion, we get using basic Brownian excursion theory, the
definition of the map $\lambda_a$ and a change of variable, that
\begin{eqnarray*}
S_1 &=& E\int\limits_0^\infty d\tilde L_s^0 \int\limits_{U}\tilde H(s,\omega,u)~ n(du)\\
&=& E\int\limits_0^\infty d\tilde L_s^0 \int\limits_{U_0^+} H(s+T_a,\omega,\lambda^a(u))I_{\{T_a < \infty\}}~ n(du)\\
&=& E \int\limits_0^\infty dL^a_{s+T_a} \int\limits_{U}
H(s+T_a,\omega,u) 1_{\{a\}} (W_{s+T_a})~n_a (du)
\end{eqnarray*}
where we have used the fact that $\tilde L_s^0= L^a_{s+T_a}$, and
that the latter process is, almost surely, supported on the set
$\{s:W_{s+T_a}=a\}$ to obtain the last equality. Hence
\[
S_1 =E\int\limits_0^\infty dL^a_s \int\limits_{U} H(s,\omega,u)
1_{\{a\}} (W_s) n_a (du).
\]
Similarly,
\[
S_2 =E\int\limits_0^\infty dL_s^b \int\limits_{U} H(s,\omega,u)
1_{\{b\}} (W_s) n_b (du).
\]
Hence
\begin{eqnarray*}
S_1 +S_2 &=& E\int\limits_0^\infty dL_s^a \int\limits_U H(s,\omega,u)~n(W_s,du)\\
&& + E\int\limits_0^\infty dL_s^b \int\limits_U H(s,\omega,u)~ n(W_s,du) \\
&=& E\int\limits_0^\infty d(L^a_s +L^b_s) \int\limits_U H(s,\omega,u)~ n(W_s,du)\\
&=& E\int\limits_0^\infty dL_s \int\limits_U H(s,\omega,u)~ n(W_s,du)\\
&=& E\int\limits_0^\infty ds \int\limits_U H(\tau_s,\omega,u)~
n(W_{\tau_s},du).
\end{eqnarray*} This completes the proof of the Theorem.\hfill{$\Box$}\\
\begin{Remark}\rm{ We can arrive at the above result , using the results
in \cite{M}. Consider the closed homogenous set $M = \{t : W_t = a~
or~ b\}$ and the corresponding exit system
$\{(\tilde{L},\tilde{n}(x,.)) ; x \in {\mathbb R} \}$. Then
$\tilde{L}_{\cdot} = L_{\cdot} = L^a_{\cdot} + L^b_{\cdot}$ and
$\tilde{n}(x,.)$ is a measure on $\tilde{U}:= U(S)\cup U(S_1)\cup
U(S_2)$ where $S:= (a,b), S_1:= (-\infty ,a), S_2 := (b,\infty)$ and
$\tilde{n}(a,.)|_{U(S)} = n(a,.) , \tilde{n}(b,.)|_{U(S)} =
n(b,.)$.}
\end{Remark}
\section{Excursions Stradling a fixed time :} We next look at
excursions into $S$, straddling a given time $t
> 0$. We work with the filtration generated by the Brownian motion $(W_t)$.
In other words, ${\cal F}_t$ is the same as ${\cal F}^W_t := \sigma
\{W_s, s \leq t\}$ augmented by all $P$ null sets. Recall that
$\sigma_t := \sup\{s \leq t : W_s \notin S \}$. Let ${\cal
F}_{\sigma_t} := \sigma \{ H (\sigma_t): \rm {H(t,\omega)~an~ {\cal
F}_t ~optional~process}\}$. Recall that $i_t(\omega)$ denotes the
portion of the excursion straddling  $t $ , from $t$  upto its
lifetime $d_t(\omega)$ and that for  $u \in U,$ $R(u) :=  \inf \{s
>0 : u_s \notin S \}$ denotes the lifetime of the excursion $u.$
Define for $s > 0$ and $ F : U_{0,\delta} \rightarrow [0,\infty]$
measurable, with $F(\delta)= 0,$ \bea q(s,F) := \frac{1}{n\{R >
s\}}\int\limits_{\{R > s\}}F(u)\,n(du). \eea
 Note that $n\{R > s\} > 0 $ (see Remark 2.1).  We then have the
following proposition from  Chap.XII, Propn.3.3, \cite{RY}:
\begin{Proposition} Let $S = (-\infty,0)\bigcup(0,\infty)$ and $\sigma_t$
the associated last entrance time before t for S. Then for every $t
> 0,$ \bea E[F(i_{\sigma_t})|{\cal F}_{\sigma_t}] =
q(t-\sigma_t,F)~a.s. \eea
\end{Proposition}

We now wish to generalise this proposition to the case of excursions
into $(a,b)$ straddling $t >0$. For $s>0, x \in \mathbb{R}$, let
\bea q(x,s,F)=\frac{1}{n(x,\{R
>s\})} \int\limits_{\{R >s\}} F(u) n(x,du). \eea

\begin{Theorem} Let $S=(a,b)$ and $U_{\delta},\sigma_t,{\cal F}_{\sigma_t}$
be associated with $S$ as above. Let $F : U_{\delta} \rightarrow
[0,\infty]$ be measurable
 with $F(\delta) = 0.$
Then for every $t >0,$ we have
 \bea E\left[F \left(i_{\sigma_t}(\omega)\right) \mid {\cal F}_{\sigma_t}\right] =
 q \left(W_{\sigma_t}, t-\sigma_t, F\right) \eea almost surely on $(\sigma_t < t)$.
\end{Theorem}

The proof of the Theorem 3.2 depends on an extension of Proposition
3.1 which we now formulate . Let $\tilde (W_s) \equiv ((W_{s+T_a}
-W_{T_a})I_{\{T_a < \infty\}})$ be the standard Brownian motion
introduced in the proof of Theorem 2.3 and $\tilde {\sigma}_t =\sup
\{s \leq t: \tilde W_s=0\}.$ Let \Bea \tilde{T}_{b-a}^t &:=& \inf\{s
> 0 : \tilde{W}_{s+ \tilde{\sigma}_t} = b-a\}~~ {\rm on}~~
\{\tilde{\sigma}_t < t \} \\ &=& \infty~~ {\rm on}~~
\{\tilde{\sigma}_t = t\}. \Eea Similarly let $(\widehat W_s) \equiv
((W_{s+T_b} -W_{T_b})I_{\{T_b < \infty\}})$, $\widehat{\sigma}_t :=
\sup \{s \leq t : \widehat{W}_{s} =0\}$ and ~\Bea
\widehat{T}_{-(b-a)}^t &:=& \inf\{s > 0 :
\widehat{W}_{s+\widehat{\sigma}_t} = -(b-a)\}~~{\rm
on}~~\{\widehat{\sigma}_t < t\} \\ & =& \infty ~~{\rm on}~~
\{\widehat{\sigma}_t = t\}.\Eea In what follows, we will abuse
notation to refer to $i_{\tilde{\sigma}_t}$ as the excursion of
$(\tilde{W}_s)$ from zero, starting at time $\tilde{\sigma}_t < t$
and a similar reference to $i_{\widehat{\sigma}_t}$ will mean the
excursion of $(\widehat{W}_s)$ from zero starting at time
$\widehat{\sigma}_t < t$. For $ u \in U_{0}$, again by abusing
notation we will denote by $T_{b-a}(u)$ the hitting time of level
$b-a$ by the excursion $u$ with a similar convention  for
$T_{-(b-a)}(u)$. The following Proposition relates the excursions
(stradling t) of $(\widehat{W}_s) $ and $(\tilde{W}_s)$ below and
above zero with state spaces $U_0^+ ,U_0^-$ respectively to the
excursions (stradling t) of $(W_s)$ into $S = (a,b)$ with state
space  $U = U(S) = U(a,b)$.

\begin{Proposition} Let $F : U_{\delta} \rightarrow [0,\infty]$ be measurable
with $F(\delta) = 0$. Let $ t > 0$ be fixed.
$$a)~~~~ E[1_{(0,t)} (\tilde{\sigma}_t) 1_{(t-\tilde{\sigma}_t,
\infty)}(T_{b-a}( i_{\tilde{\sigma}_t}))F\circ
\lambda^a(i_{\tilde{\sigma}_t})|{{\cal F}}_{\tilde{\sigma}_t}] = q(
W_{\sigma_t},t-\sigma_t,F)
$$ ~a.s. on the set $(\tilde{\sigma}_t< t,
\tilde{T}_{b-a}^t > t- \tilde{\sigma}_t)$ and \\ $$b)~~~~
E[1_{(0,t)}(\widehat{\sigma}_t )1_{(t-\widehat{\sigma}_t,
\infty)}({T}_{-(b-a)}(i_{\widehat{\sigma}_t}))F\circ
\lambda^b(i_{\widehat{\sigma}_t})|{{\cal F}}_{\widehat{\sigma}_t}] =
q(W_{\sigma_t},t-\sigma_t,F)
$$~a.s. on the set $(\widehat{\sigma}_t <t,
\widehat{T}_{-(b-a)}^t
> t - \widehat{\sigma}_t).$
\end{Proposition}
{\bf Proof :}Let $\alpha (s,\omega)$ be $({\cal F}_t)$-optional. To
prove a) we need to show that \Bea && E[ \alpha (\tilde{\sigma}_t)
\, 1_{(0,t)} (\tilde{\sigma}_t) 1_{(t-\tilde{\sigma}_t,
\infty)}(T_{b-a}( i_{\tilde{\sigma}_t})) \, F(\lambda^a (
i_{\tilde{\sigma}_t}))]
\\  &=& E[\alpha (\tilde{\sigma}_t) \, q(W_{\sigma_t},
t-{\sigma}_t, F)\,1_{(0,t)} (\tilde{\sigma}_t)
1_{(t-\tilde{\sigma}_t, \infty)}(\tilde{T}_{b-a}^t)]\Eea For
$(s,\omega,u) \in [0,\infty) \times \Omega \times U_{0,\delta}$
define \Bea H(s,\omega,u) &:=& \alpha (s) ~F \circ \lambda^a
(u)~I_{(0,t)} (s)~ I_{\{R > t-s\}} (u) ~I_{\{T_{b-a} > t-s\}} (u).
\\\Eea
Let $\tilde G (\omega) \subset [0,\infty)$ be the left end points of
excursion intervals of $(\tilde W_t)$ from zero. Let $L_t :=
\tilde{L}_t^0,$ the local time of zero of $(\tilde {W}_s)$. Recall
that $(\tau_t)$ is the right continuous inverse of $(L_t)$. Then we
may write as in Proposition 3.3, Chapter XII, \cite{RY},\Bea
&&E\left[ \alpha (\tilde \sigma_t) I_{(0,t)} (\tilde \sigma_t)
I_{(t-\tilde \sigma_t,\infty)} \left( T_{b-a} (i_{\tilde
\sigma_t})\right) F\circ \lambda^a (i_{\tilde
\sigma_t})\right]\\
&=& E \sum\limits_{s\in \tilde G} H(s,\omega,i_s (\omega))\\
&=& E \int\limits_0^\infty ds \int\limits_{U_{0,\delta}} H(\tau_s, \omega, u)~ n (du)\\
&=& E\int\limits_0^\infty ds ~I_{(0,t)} (\tau_s) ~\alpha (\tau_s)
\int\limits_{\{R > t- \tau_s\} \cap \{T_{b-a} > t- \tau_s\}}F\circ \lambda^a (u)~n(du)\\
&=& E \int\limits_0^\infty ds ~I_{(0,t)} (\tau_s) ~q  (a,t-\tau_s,F) n^a (R > t-\tau _s)\\
&=& E \sum\limits_{s\in \tilde G} \alpha (s) ~q(a,t-s,F)~I_{(0,t)} (s) G(s,i_s)\\
&=& E~\alpha (\tilde \sigma_t) ~I_{(0,t)} (\tilde \sigma_t)
~q(a,t-\tilde \sigma_t,F)~G(\tilde \sigma_t, i_{\tilde \sigma_t})\\
\Eea where for $(s,u) \in [0,\infty) \times U_{0,\delta}$ we define
\Bea G(s,u) := I_{U^+_0 \cap \{R > t-s, T_{b-a} > t-s\}}(u); \Eea
and where in the 4th equality above we have used the fact that for
$0 < s <t$, \Bea (\lambda^a)^{-1} \{u \in U :R(u) > t-s\} =U_0^+
\cap \{R > t-s, T_{b-a} > t-s\}. \Eea Finally we note that when
$G(\tilde \sigma_t, i_{\tilde \sigma_t})=1$, $\tilde \sigma_t
=\sigma_t$ and $W_{\sigma_t} =a$. This completes the proof of a).
The proof
of b) is similar. \hfill{$\Box$}\\

 {\bf Proof of Theorem 3.2:} Let $\alpha (s,\omega)$ be $({\cal F}_t)$-optional
 and $t > 0$ be given.
We need to show \Bea E \alpha \left(\sigma_t \right)
1_{\left(\sigma_t < t \right)} F\left(i_{\sigma_t}\right) = E\alpha
\left(\sigma_t \right) 1_{\left(\sigma_t <t\right)}
q\left(W_{\sigma_t}, t-\sigma_t, F\right). \Eea  Recalling the
notation $i_{\tilde{\sigma}_t}$ and $i_{\widehat{\sigma}_t}$
introduced before the statement of Proposition(3.3), we note that
\Bea (\sigma_t <t) &=&(\tilde{\sigma}_t< t, \tilde{T}_{b-a}^t > t-
\tilde{\sigma}_t)\cup (\widehat{\sigma}_t <t, \widehat{T}_{-(b-a)}^t
> t - \widehat{\sigma}_t)\\ &=&  (\tilde{\sigma}_t< t,
{T}_{b-a}(i_{\tilde{\sigma}_t}) > t- \tilde{\sigma}_t)\cup
(\widehat{\sigma}_t <t, {T}_{-(b-a)}(i_{\widehat{\sigma}_t})
> t - \widehat{\sigma}_t)\Eea where the sets in the right hand side
 are disjoint. Further we note that
 $\sigma_t \,
 =\tilde{\sigma}_t \,$  on the set $(\tilde{\sigma}_t< t,
\tilde{T}_{b-a}^t > t- \tilde{\sigma}_t)$ and $\sigma_t =
\widehat{\sigma}_t$ on the set $(\widehat{\sigma}_t <t,
\widehat{T}_{-(b-a)}^t
> t - \widehat{\sigma}_t).$ We then have \Bea  E \alpha \left(\sigma_t \right)
1_{\left(\sigma_t < t \right)} F\left(i_{\sigma_t}\right) &=& E[
\alpha (\tilde{\sigma}_t) \, 1_{(0,t)} (\tilde{\sigma}_t)
1_{(t-\tilde{\sigma}_t, \infty)}(T_{b-a}( i_{\tilde{\sigma}_t})) \,
F(\lambda^a ( i_{\tilde{\sigma}_t}))] \\ &+& E[ \alpha
(\widehat{\sigma}_t) \, 1_{(0,t)} (\widehat{\sigma}_t)
1_{(t-\widehat{\sigma}_t, \infty)}(T_{-(b-a)}(
i_{\widehat{\sigma}_t})) \, F(\lambda^b ( i_{\widehat{\sigma}_t}))]\\
&=& E[\alpha (\tilde{\sigma}_t) \, q(W_{\sigma_t}, t-{\sigma}_t,
F)\,1_{(0,t)} (\tilde{\sigma}_t) 1_{(t-\tilde{\sigma}_t,
\infty)}(\tilde{T}_{b-a}^t)] \\ &+& E[\alpha (\widehat{\sigma}_t) \,
q(W_{\sigma_t}, t-{\sigma}_t, F)\,1_{(0,t)} (\widehat{\sigma}_t)
1_{(t-\widehat{\sigma}_t, \infty)}(\widehat{T}_{b-a}^t)] \\ &=&
E\alpha \left(\sigma_t \right) 1_{\left(\sigma_t <t\right)}
q\left(W_{\sigma_t}, t-\sigma_t, F\right), \Eea where to obtain the
second equality we have used the result of Proposition(3.3). This
completes the proof of Theorem(3.2). \hfill {$\Box$}

\section{Asymptotic Distribution of Excursions \\ Stradling a fixed time :}

Let $C := C([0,\infty),\bar{S}) \cup \{\delta\}$ with the
$\sigma$-field ${\cal C}$ generated by the Borel sigma field of
$C([0,\infty),\bar{S})$ and the singleton $\{\delta\}$. We now
consider only the case $S = (a,b).$ All excursion related objects
are considered with respect to this $S.$ We consider the excursions
$i_{\sigma_t}$ as a stochastic process with values in $ C $ and
accordingly use a new notation. We define the $C$ valued stochastic
process $(\zeta_t)$, measurable in $(t,\omega)$ as follows: \Bea
\zeta_t(.) := W_{\sigma_t + \cdot\wedge (d_t -\sigma_t)} \,
I_{(\sigma_t <t)} +\delta \, I_{(\sigma_t =t)} \Eea where on
$(\sigma_t <t)$, we note that the function $s \rightarrow
W_{\sigma_t +s \wedge (d_t -\sigma_t)}$ belongs to  $C$. Let $E:=
\{a,b\}\times[0,\infty)\times C$ and ${\cal E}:=
\{a,b,\phi,\{a,b\}\}\times{\cal B}[0,\infty)\times{\cal C}$ the
product sigma field on $E$. For $A \in {\cal E}, x=a,b, ~{~\rm and~}
s>0$, let $A(x,s):= \{\omega: (x,s,\omega) \in A\}$. For $A \in
{\cal E}, s
> 0, x \in \R$, we recall from eqn.(3.6) the kernel $$ q(x,s,A(x,s)) : = \frac
{n(x,A(x,s)\bigcap \{R > s\})}{n(x,\{R > s\})}.$$ We then define the
probability measure $P^0$ on $(E,{\cal E})$ as follows : \bea P^0
(A) : = \int\limits_0^\infty
\left(\frac{q(a,s,A(a,s))+q(b,s,A(b,s))}{2}\right)\, dF(s) \eea
where $F(\cdot)$ is a distribution function on $[0,\infty)$ defined
as follows : Let $P_x , x \in \R ,$ denote the distribution of $(W_s
+ x)$ on $C([0,\infty),\R).$ Then, \Bea F(s) := \frac{1}{b-a}
\int\limits_a^b (1-\psi (x,s))\, dx \Eea where \Bea \psi(x,s) := P_x
\left( a < \inf\limits_{0 \leq r \leq s} W_r < \sup\limits_{0 \leq r
\leq s} W_r < b\right). \Eea We note that $F$ is the asymptotic
distribution of $t-\sigma_t$, conditional on \\$\{W_t \in (a,b)\}$
as $t \rightarrow \infty$ (see \cite{AR0},Thm.(4.2)). Further it is
clear that $P^0(A) =E\,q(X,Y,A(X,Y))$ where $X,Y$ are independent,
$Y \sim F$ and $P(X=a)=P(X=b)=\frac{1}{2}$.

We then have the following theorem.
\begin{Theorem} Let $(\zeta_t), P^0,X,Y
$ be as above. Then, conditional on $(\sigma_t < t),
(W_{\sigma_t},t-\sigma_t,\zeta_t)$ converges weakly to $P^0$ on $
(E,{\cal E})$ as $t \rightarrow \infty$.
\end{Theorem}

{\bf Proof:} Let $f: E \rightarrow \mathbb{R}$ be a bounded and
continuous function. It suffices to show that $$
\lim\limits_{t\rightarrow \infty} \,
E[f(W_{\sigma_t},t-\sigma_t,\zeta_t)\mid \sigma_t <t] =
\int\limits_{E} f\, dP^0.$$ We have  \Bea
E[f(W_{\sigma_t},t-\sigma_t,\zeta_t) \mid \sigma_t < t] &=&
\frac{E\,I_{(\sigma_t <t)}
f(W_{\sigma_t},t-\sigma_t,\zeta_t)}{P(\sigma_t < t)}\\
&=& \frac{E\left[ I_{(\sigma_t < t)}
E[f(W_{\sigma_t},t-\sigma_t,\zeta_t) \mid {\cal
F}_{\sigma_t}]\right]}{P(\sigma_t <t)}. \Eea From Theorem 3.2, we
have \Bea E[f(W_{\sigma_t},t-\sigma_t,\zeta_t) \mid {\cal
F}_{\sigma_t}]=q(W_{\sigma_t},
t-\sigma_t,f(W_{\sigma_t},t-\sigma_t,\cdot)) \Eea almost surely on
$\{\sigma_t < t\}$. On the other hand we know from the results in
\cite{AR} that $(W_{\sigma_t} , t-\sigma_t)$ converges weakly to
$(X,Y)$ conditional on $\{\sigma_t <t\}$. Using Remark 2.1 it can be
shown that $q(x,s,f(x,s,\cdot))$ is a bounded continuous function on
$\{a,b\} \times (0,\infty)$. The result follows. \hfill{$\Box$}

\begin{Remark} \rm {The limiting distribution $P^0$ can be described in
terms of the measure $Q_{ex}$ introduced in \cite{P}. We recall that
the measure $Q_{ex}$ was introduced on the space of excursions as an
`equilibrium measure' or more specifically as the `Palm measure'
corresponding to a stationary point process. It is natural to
interpret our results in Sec.4 as $ t \rightarrow \infty $ in terms
of this equilibrium measure.
   Let $ M , \{(\tilde{L},\tilde{n}(x,.)) ; x = a~ or~ b \}$ be as in Remark
2.4, with $\tilde{L}_{\cdot} = L_{\cdot} = L^a_{\cdot} +
L^b_{\cdot}$. Then it was shown in \cite{P} that $Q_{ex}$  was given
as $ Q_{ex}(.) = \int \alpha(dx)\tilde{n}(x,.)$ where $\alpha(A) :=
\int dx E^x \int_0^1I_A(W_s)dL_s$. It follows from the occupation
density formula that $\alpha = \delta_a + \delta_b$ and consequently
from Remark 2.4 that, $Q_{ex}|_{U(S)} = n(a,.)+ n(b,.)$. We than
have for $x = a~ or ~b , s > 0, A \in {\cal U}(S)$
$$ q(x,s,A) : = \frac
{Q_{ex}(A\bigcap \{R(u) > s, u(0)=x\})}{Q_{ex}(\{R(u) >
s,u(0)=x\}}.$$}
\end{Remark}
{\bf An application :} For $u > 0, y > 0,$ consider the
probabilities \Bea \phi(t) := P_a\{ 0 < t - \sigma_t < u, 0 < W_t -
W_{\sigma_t} < y\} \Eea where for each $a \in {\mathbb
R},~P_a\{\omega:W_0 = a\}= 1$. In \cite{AR0},Theorem (5.1), it was
shown that $\phi(t)$ satisfies a renewal equation. With the renewal
theorem in mind, a natural question is to find the limit of
$\phi(t)$ as $t \rightarrow \infty$. However, it is easy to see that
the probabilities that define $\phi(t)$ converge to zero since the
events in question are contained in $\{W_t \in (a,b)\}$. On the
other hand, since these events can be expressed as functionals of
the excursion straddling $t$, we can apply the previous theorem to
compute the limit of the conditional probability
$$ P_a\{\omega : 0 < t - \sigma_t < u, 0 < W_t - W_{\sigma_t} < y| W_t
\in (a,b)\}$$ as $t \rightarrow \infty$. Let
$$A := \{(x,s,\omega): 0 < \omega_s - \omega_{0} < y, 0<s<u, x=a~
{\rm or}~ b\}.$$ Then by Thm.(4.1),  \Bea && P_a\{0<t - \sigma_t <
u, 0 < W_t -
W_{\sigma_t} < y| W_t \in (a,b)\} \\
&=& P_a \{(W_{\sigma_t},t-\sigma_t,\zeta_t) \in A| \sigma_t < t\}\\
&\rightarrow & P^0(A)\Eea provided $P^0(\partial A) = 0$.  Note that
\Bea (\partial A) &=& (\bar{A}-A^\circ) =  \{(x,s,\omega): s=u,0
\leq \omega_u - \omega_0 \leq y, x= a{\rm~or~}b\}\\ && \bigcup
\{(x,s,\omega): 0 < s < u, \omega_s -\omega_0 = 0~{\rm or~}y,x=
a{\rm~or~}b\}\\ &=:& A_1\bigcup A_2 \Eea Clearly, from the
definition of the measure $P^0$ we have,
$$P^0(A_1)= P^0\{(a,s,\omega): s=u,0 \leq \omega_u - \omega_0 \leq y,x= a{\rm~or~}b\} = 0.$$ As for
the second set $A_2$ in the union, from the definition of $P^0$
(eqn.(4.8)) it suffices to show that for every $0 < s <u, x= a,b,
q(x,s,(A_2)(x,s))=0$. Again, from the definition of $q(x,s,.)$
(eqn.(3.6)) it suffices to show that for each $x=a,b$ and $0 < s <
u, n(x,A_2(x,s)\bigcap \{R > s\})=0$. From the definition of the
kernels $n(x,\cdot)$, we have for $0 < s < u, x= a$, \Bea
n(a,\{\omega&:& \omega_s -\omega_0 = 0~{\rm or~}y\}\bigcap \{R
> s\}) \\
&=& n^+\circ \lambda_a^{-1}(\{\omega :\omega_s -\omega_0 = 0~{\rm
or~}y\}\bigcap \{R > s\}) \\ &=&  n^+(\{\omega :\omega(s\wedge
T_{b-a}) = y\}\bigcap \{R\circ \lambda_a
> s\}) \\ &=&  n^+(\{\omega :\omega(s\wedge
T_{b-a}) = y\}\bigcap \{R_1\wedge T_{b-a}
> s\}) \\ &=& n^+(\{\omega :\omega(s) = y\}\bigcap \{R_1\wedge T_{b-a}
> s\}) \\ &=& 0.  \Eea
where $R_1(\omega)$and $T_{b-a}(\omega)$ are respectively the life
time and the hitting time of $b-a$ of the excursion $\omega$
starting at $0$ and the last equality follows from the absolute
continuity of the map $n\circ \omega_s^{-1}, s > 0$. This proves
that $P^0(\partial A)= 0$.

{\bf Acknowledgement :} The author would like  to thank Jean Bertoin
for pointing out reference \cite{P} .

\end{document}